\documentclass[11pt]{article}

\usepackage{graphicx}
\usepackage[english]{babel}
\usepackage{authblk}
\usepackage{subcaption}  % for subfigures
\usepackage{caption} 
\usepackage{multirow}

\usepackage[table]{xcolor}
\usepackage[colorlinks = true,
            linkcolor = blue,
            urlcolor  = blue,
            citecolor = blue,
            anchorcolor = blue]{hyperref}

\usepackage[margin=2.5cm]{geometry}
 \usepackage{amsmath, amssymb}
\usepackage{amsthm}

\usepackage{aliascnt}

\newaliascnt{remark}{theorem}
\newtheorem{remark}[remark]{Remark}
\aliascntresetthe{remark}

\newaliascnt{example}{theorem}

\aliascntresetthe{example}

\newcommand{\R}{\mathbb{R}}
\newcommand{\N}{\mathbb{N}}

\begin{document}

\title{Estimating parameters of the diffusion model via asymptotic expansions}

\author[1]{K. Kalimeris\thanks{kkalimeris@Academyofathens.gr}}
\affil[1]{Mathematics Research Center, Academy of Athens, Greece}
\author[2]{L. Mindrinos\thanks{leonidas.mindrinos@aua.gr}}
\affil[2]{Department of Natural Resources Development and Agricultural Engineering, Agricultural University of Athens, Greece}

\maketitle

\begin{abstract}
A broad class of inverse problems deals with determining certain parameters, from measurement data, in models which are associated to certain partial differential equations. In this work we focus on the heat equation on a finite interval and we determine the dimensionless diffusion parameter from a single measurement. Our results extend to estimating additional parameters of the initial-boundary value problem, such as the length of the interval and/or the time required for the solution to achieve a specific state. Our approach relies on the asymptotic solution of an integral equation: The formulation of this integral equation is based on the solution of the direct problem via the Fokas method; the solution of this equation is achieved through the asymptotic evaluation of the associated integrals which yield an effective approximate solution, supported by numerical verifications. We apply these approximations to well-established problems in soil science and we compare our results with existing ones, displaying clear improvement.
\end{abstract}

\section{Introduction and Problem formulation}

The inverse problem of parameter estimation concerns the recovery of unknown parameters or functions from exact or noisy measurement data. Such problems frequently arise in physical and engineering applications, where the governing processes are modelled by initial boundary value problems (IBVPs) or integral equations. Because of their non-linearity and ill-posedness, these problems exhibit strong sensitivity to data perturbations, often leading to unstable or non-unique solutions. In the case of clean data, exact reconstruction of the parameters may be feasible; however, when noise is present, the problem is typically reformulated as an optimization problem involving the associated forward operator. In the present work, we are interested in introducing a methodology to deal with a family of problems with noiseless data; similar problems can be found in \cite{dantas, fadale, mohebbi}. Addressing the ill-posedness, characterized by the lack of stability, will be considered in future work involving regularization techniques which are used to ensure stable and physically meaningful parameter estimates; as indicative literature we consider the classic monographs \cite{kunisch, engl, barbara}, and the references therein.

In the current manuscript we employ results from the Fokas method \cite{fokas97, fokas08} as the starting point for the formulation of the inverse problem. In 1967 in the seminal paper \cite{IST67} the method of inverse scattering transform (IST) was introduced for the solution of initial value problems (IVPs) for integrable partial differential equations (PDEs). In 1997, Fokas introduced the  unified transform method, which now is known as the Fokas method,  for the solution of IBVPs for integrable PDEs, providing a positive answer to that long-standing extension of IST to problems with prescribed boundary conditions; later, in \cite{fokas02b}, it was first realised that novel forms of solutions can be provided by this method for linear PDEs, too. In the current manuscript, we employ results from \cite{arg24} and \cite{KalMin25b}, where the authors provide the integral representation of the solution to two IBVPs, which are associated to physical problems in agriculture for bounded domains; this is also to be contrasted to the problem formulation via the traditional methods for treating these problems, where a series representation of the solution is involved. 

For a diffusion model, the estimation of the physical parameters relies on prescribing some known values on instances of the solution of the associated direct problems. Applying the above-mentioned methodology in this procedure, gives rise to an integral equation. Once the integral equation is constructed, we employ techniques from the broad area of  asymptotic expansions to properly approximate its solution, see for example \cite{erdelyi56}. This include asymptotic estimate of the contribution of integrals \cite{bleistein75, AF03}, as well as other approaches involving  consistent balancing of asymptotic expansions, which can be found in \cite{lag13, AF03, miller06}. To our knowledge, this methodology has never appeared in the literature, since the relevant works are not formulating an integral equation: the origin of this difference comes from the fact that the solution of the direct problem, via the traditional methods, is given in terms of (Fourier-type) series. The fact that IBVPs of bounded domains admit discrete spectrum (resp. series representation of the solution) is a broadly accepted misconception, but the development of the Fokas method reveals the image which is closer to the truth, namely that IBVPs of both bounded and unbounded domains admit continuous spectra, (resp. integral representation of the solution), see \cite{CFK25}; for a rigorous analysis of this equivalence we refer to the recent work \cite{OMKarx}; this concept is also true for the case of space and time dependent parameters of the problem \cite{FD23, DF25, KOarx}. This uniformity was initially discussed in \cite{FP05}, where it was also observed that there are linear and bounded problems that admit \textit{only} continuous spectra.  

We inspect the performance of the associated asymptotic expansions in the estimation of the model parameters in \textit{three} ways. First, in the mathematical approach the order of magnitude of the remainder of the truncated asymptotic expansion is considered. Second, we provide illustrations of the numerical results  as a verification which involves the relative error of the estimated value versus the real value of the relevant parameter. Third, we investigate the performance of our algorithm for the estimation of physical parameters in specific problems of agricultural engineering. In soil science, numerous mathematically interesting inverse problems can be formulated as the solution of an integral equation, due to their correspondence to IBVPs which model the associated direct problem. The central objective in such problems is to infer soil properties or other unknown parameters from observed hydrological states.

The first two inspections of the above paragraph are presented in \autoref{sec_approx}, where several techniques from the asymptotic analysis are employed in order to cover the cases of estimating both large and small values of the associated parameters. In \autoref{sec_main}, a proper matching of asymptotics is provided, so that a combination of the results of  \autoref{sec_approx} would yield optimal results, balancing accuracy of estimates with easy-to-apply formulas for the parameters of the model. The third inspection is provided in \autoref{sec_num} of this study, where we focus on two distinct inverse problems and develop a unified framework for the recovery of the parameters of interest. Although these problems refer to different physical scenarios and originate from different  PDEs, under certain assumptions they can both be reduced to an IBVP for the heat equation (namely the diffusion model), but with different initial and boundary conditions. Thus, the inverse problems can be expressed in terms of an integral equation that exhibits the same fundamental properties which admits a coherent reconstruction method. 

In the rest of the current section we show the reduction of those two physical problems to the mathematical formulation of the typical integral equation \eqref{eq_final_inverse}, in view of the definition \eqref{eq_integral}, which is the constitutional equation that we will asymptotically solve in \autoref{sec_approx}.

\subsection{Estimating drain spacing from water table observations}\label{sec_drain}

In agricultural engineering, determining the optimal spacing of subsurface drains is critical for efficient water management. The inverse problem in this context involves estimating the drain spacing required to achieve a target water table depth based on desired water table elevations over time.

In \cite{KalMin25b} the authors derived an analytical solution to the forward problem, allowing computation of the water table height from known soil properties, under the assumption of horizontal flow only. The corresponding IBVP is formulated as follows:
\begin{subequations}\label{bvp_heat}
\begin{alignat}{3}
 \frac{\partial h}{\partial t} &= \mathcal{A}\frac{\partial^2 h}{\partial x^2},  \quad && 0<x <L, \, t>0,  \label{bvp1h}\\
h (x,0) &= d+ h_0, \quad && 0 <x <L, \label{bvp2h}\\ 
h (0,t)  &= d,  \quad &&t>0, \label{bvp3h}\\
 \frac{\partial h}{\partial x} (L,t) &= 0,  \quad &&t>0, \label{bvp4h}
\end{alignat}
\end{subequations}
where $h(x,t)$ denotes the water table height, $\mathcal{A}$ is a positive constant, and $d$ represents the known elevation of the drain axis above the impervious layer. In this study, the initial water table, $h_0,$ is assumed to be constant, representing a scenario in which the groundwater is initially uniform.

Using the Fokas method \cite{fokas97, fokas08}  the solution of \eqref{bvp_heat} is given by \cite[Eq. (20)]{KalMin25b} 
\begin{equation}\label{eq:20}
h(x,t) = d + h_0 -i \frac{h_0}{\pi} \int_{C} e^{-\mathcal{A} \lambda^{2} t} \frac{\cos[\lambda(x - L)]}{\lambda \cos(\lambda L)}  d\lambda,
\end{equation}
with $C$ being a suitable curve in the complex $\lambda$-plane, under the instructions given in \cite{KalMin25b}; here we choose the hyperbola parametrised as $\lambda(r)=2 \sinh r \,+ \,i \,\cosh r, \ r\in\R.$

Let $h_0, \ d$ be fixed. The inverse problem in mathematical terms reads: given $T$ and $\mathcal{A}$, find $L$, such that
\begin{equation}\label{eq:condition}
h(L,T) = H \quad \text{(known)}.
\end{equation}
We rewrite \eqref{eq:20} with condition \eqref{eq:condition} as follows:
\begin{equation}\label{eq:rewritten}
\frac{i}{\pi} \int_{C} e^{-\mathcal{A} \lambda^{2} T} \frac{1}{\lambda \cos(\lambda L)}  d\lambda = \frac{d - H}{h_0} + 1.
\end{equation}
Let the change of variables \( \lambda L = k \), then
\begin{equation}\label{eq:dimensionless}
\frac{i}{\pi} \int_{C} e^{-\frac{\mathcal{A} T}{L^2} k^{2}} \frac{1}{k \cos k}  dk = \frac{d - H}{h_0} + 1,
\end{equation}
where $C$ remained invariant under this change of variables, under the associated contour deformation.

We define
\begin{equation} \label{def:c+a}
c_1 := \frac{d - H}{h_0} + 1, \quad \text{ and } \quad a_1 := \frac{\mathcal{A} T}{L^2}.
\end{equation}
Being consistent with the underlying physical model, the admissible parameters are restricted to
 \begin{equation} \label{restrictions-on-H}
d < H < d+h_0 \quad \Leftrightarrow \quad 0 < \frac{d - H}{h_0} + 1 < 1,  
\end{equation}
since, after some time $T>0$, the height $H$ will have decreased below its initial value $d + h_0$, yet remain above the minimum $d$ observed at the left boundary. Thus,  \( c_1 \in (0,1) \), 
 defined by \( c_1 := c_1 \left( \frac{\mathcal{A} T}{L^2} \right) = c_1 (a_1) \), with \( c_1 \rightarrow 0\) as $a_1 \rightarrow 0$, and \( c_1 \rightarrow 1\) as $a_1 \rightarrow \infty,$ verifying  the maximum principle for the heat equation.

We define
\begin{equation}\label{eq_integral}
I (a) :=  \frac{i}{\pi} \int_{C} e^{-a k^{2}} \frac{1}{k \cos k}  dk
\end{equation}
and the constitutional integral equation reads
\begin{equation}\label{eq:Ia}
I(a_1 ) = c_1,
\end{equation}
for \( a_1 > 0 \), and \( c_1 \in (0,1) \).

\subsection{Determining soil properties from moisture data}\label{sec_infil}

The inverse problem involves estimating soil hydraulic parameters--such as saturated hydraulic conductivity ($K$) and diffusivity ($D$)--from measurements of soil moisture at various depths.

We focus on the case of vertical infiltration, corresponding to a one-dimensional problem, in a bounded soil profile subject to flooding (i.e., Dirichlet boundary conditions). The soil is assumed to be homogeneous with constant diffusivity, $D = D_0 > 0$, and gravity effects are neglected, $K = 0$. Under these assumptions, the mathematical model describing the problem is given by:
\begin{subequations}\label{bc1u}
\begin{alignat}{3}
\frac{\partial \theta}{\partial t}  &= D_0 \frac{\partial^2 \theta}{\partial x^2},  \quad && 0<x <L, \, t>0,  \\
\theta (x,0) &= \theta_0, \quad && 0 <x <L, \\ 
\theta (0,t) &= \theta_1, \quad &&t>0, \label{bc1u3}\\
\theta (L,t) &= \theta_0, \quad &&t>0, 
\end{alignat}
\end{subequations}
where $\theta (x,t)$ is the water content and $\theta_1 > \theta_0 >0$ are known constants. As shown in \cite[Section 3.1]{arg24}, the solution of \eqref{bc1u} takes the form: 
\begin{equation}\label{sol_for_infil}
\begin{aligned}
\theta (x,t) &= \theta_0 + (\theta_1 - \theta_0) \frac{2}{\pi} \int_C \frac{e^{-D_0 \lambda^2 t}}{ \lambda (e^{-i\lambda L}-e^{i \lambda L}) } \sin [\lambda(L-x) ]  d\lambda,\\
&= \theta_0 + (\theta_1 - \theta_0) \frac{i}{\pi} \int_C \frac{e^{-D_0 \lambda^2 t}}{ \lambda \sin (\lambda L) } \sin [\lambda(L-x) ]  d\lambda,
\end{aligned}
\end{equation}
where $C$ can be the same curve as in \eqref{eq:20}.

The inverse problem reads: Given a single moisture measurement at certain depth and time $T,$ find the diffusivity $D_0,$ meaning solve the non-linear integral equation
\begin{equation}\label{eq:condition2}
\theta \left(\frac{L}{2},T \right) = \Theta,
\end{equation}
for $D_0.$ Here, we assume that the data are collected at position 
$x=\frac{L}{2}$ for presentation purposes. The following analysis is valid also for other depth positions. Then, \eqref{sol_for_infil} takes the form
\begin{equation}
\frac{i}{\pi} \int_C \frac{e^{-D_0 \lambda^2 T}}{ \lambda \sin (\lambda L) } \sin \left(\frac{\lambda L}{2}\right)  d\lambda =  \frac{\Theta - \theta_0}{\theta_1 - \theta_0}.
\end{equation}
Using the identity $\sin x = 2 \sin \tfrac{x}{2} \cos \tfrac{x}{2}$ we obtain
\begin{equation}
\frac{i}{\pi} \int_C e^{-D_0 \lambda^2 T} \frac{1}{2 \lambda \cos \left(\frac{\lambda L}2 \right) }   d\lambda =  \frac{\Theta - \theta_0}{\theta_1 - \theta_0}.
\end{equation}

The change of variables \( \frac{\lambda L}2 = k \) results in
\begin{equation}\label{eq:dimensionless2}
\frac{i}{\pi} \int_{C} e^{-\frac{4 D_0 T}{L^2} k^{2}} \frac{1}{k \cos k}  dk = 2 \frac{\Theta - \theta_0}{\theta_1 - \theta_0},
\end{equation}
where $C$ remains again invariant under the change of variables. We define
\begin{equation}\label{def:c+a2}
c_2 :=  2 \frac{\Theta - \theta_0}{\theta_1 - \theta_0}, \quad \text{ and } \quad a_2 := \frac{4 D_0 T}{L^2},
\end{equation}
to get (see the 
definition \eqref{eq_integral}) the integral equation
\begin{equation}\label{eq:Ib}
I (a_2 ) = c_2,
\end{equation}
with $a_2 >0$ and $c_2 \in (0,1).$  The latter constraint is obtained by the argumentation in \autoref{sec_drain} which is based on the maximum principle of the heat equation; hence, the admissible values of $\Theta$ are such $\Theta \in \left( \theta_0, \, \frac{\theta_0+\theta_1}2\right).$

\section{Asymptotic approximation of the solution}\label{sec_approx}

%\subsection{Main Result}\label{sec_main}

We are interested in approximating the solution of 
\begin{equation}\label{eq_final_inverse}
I(a) = c, \quad a>0, \, c \in (0,1),
\end{equation}
with $I(a)$ defined in \eqref{eq_integral}, by employing several asymptotic techniques to estimate the integral $I(a),$ for large and small values of $a$, in \autoref{sec_large} and \autoref{sec_small}, respectively. This integral equation refers to both \eqref{eq:Ia} and \eqref{eq:Ib}; for the former case our results are compared to the known results of \cite{dumm1954, beers65}, displaying substantial improvement.

\subsection{Integral estimation for large values of $a$}\label{sec_large}
Integral estimation by asymptotics, for \( a \to \infty \), can be performed by the steepest descent method \cite{BH75, AF03}. In few words, deform \( C \), so that it passes through the so-called stationary point (saddle point) in the direction of the steepest descent. Here, one has to perform a modification of this approach because the stationary point $k=0$, is also a singular point of the integrand.

The first order of this approximation coincides with the residue contribution of $k=0$, which yields
\begin{equation}\label{eq:first_order}
I(a) = 1 + O(e^{-c_0 \, a}), \quad \text{for come constant } \ c_0>0,
\end{equation}
which gives exponentially small terms for \( a \to \infty \). However, it is important to find out what happens for ``small'' values of \( a>0 \), hence one should determine higher terms of this expansion. The careful book-keeping of the steepest descent method  yields the contribution produced by the singularities of the integrand, which are exponentially small, and coincide with the residue contribution of the poles $\cos k=0 \Leftrightarrow  k_n = \frac{\pi}{2} + n\pi$. In fact, one can obtain this result by deforming \( C \) to the real axis which avoids poles  $k_n$, and observe that the associated principal value integral is vanishing, by the oddness of the integrand. Then, one gets the analogue of the Fourier series \cite[Eq. (22)]{KalMin25b} and \cite{CarJae59}, namely
\begin{align}\label{eq:series2}
I(a) = 1 -\sum_{n=0}^{\infty} \frac{(-1)^n}{\left( 2n + 1 \right)\frac{\pi}{4}} \exp\left\{ -( 2n + 1)^2 \frac{\pi^2}{4} a \right\}.
\end{align}

\begin{remark}\label{remark1}
We note that the solution reported in \cite{CarJae59} corresponds to the setting of  \cite{dumm1954}, where  the problem \eqref{bvp_heat} is defined in the interval $[0,L/2]$.  When considering the inverse problem \eqref{eq:condition}, in the agricultural engineering literature \cite{dumm1954, beers65}, only the first order term of the above  Fourier series is taken into account, by imposing a threshold on the ratio of the second to the first term of the series.
% Their results corresponds to
%\\
 %if the associated parameter $\alpha \geq 0.2$
%, as the benchmark for the acceptable approximation. 
%This value of $\alpha$ corresponds to 
%$a> ....$ in our setting.
\end{remark}
Following the argumentation of the above Remark, we observe that 
%could take only the first order term, since 
the ratio of the second to the first term reads
\begin{equation}
\left| \frac{-\frac{4}{3\pi} e^{-\frac{9 \pi^2}{4}a}}{\frac{4}{\pi} e^{-\frac{ \pi^2}{4}a}} \right| = \frac{1}{3} e^{-2 \pi^2 a} \leq \frac{1}{3 \, e^8}, \quad \mbox{for} \quad a \geq \frac{4}{\pi^2}.
\end{equation}
Hence, if $a \geq \tfrac{4}{\pi^2} \approx 0.405 \, \Rightarrow \, c > 0.532$, the second term may be neglected, at the expense of an error less than $0.011\%$ of the magnitude of the first order term.% ($e^{-8} \approx 0.000335$). 

\begin{figure}[h!]
    \centering
    % First subfigure
    \begin{subfigure}[b]{0.45\textwidth}
        \centering
        \includegraphics[width=0.75\textwidth]{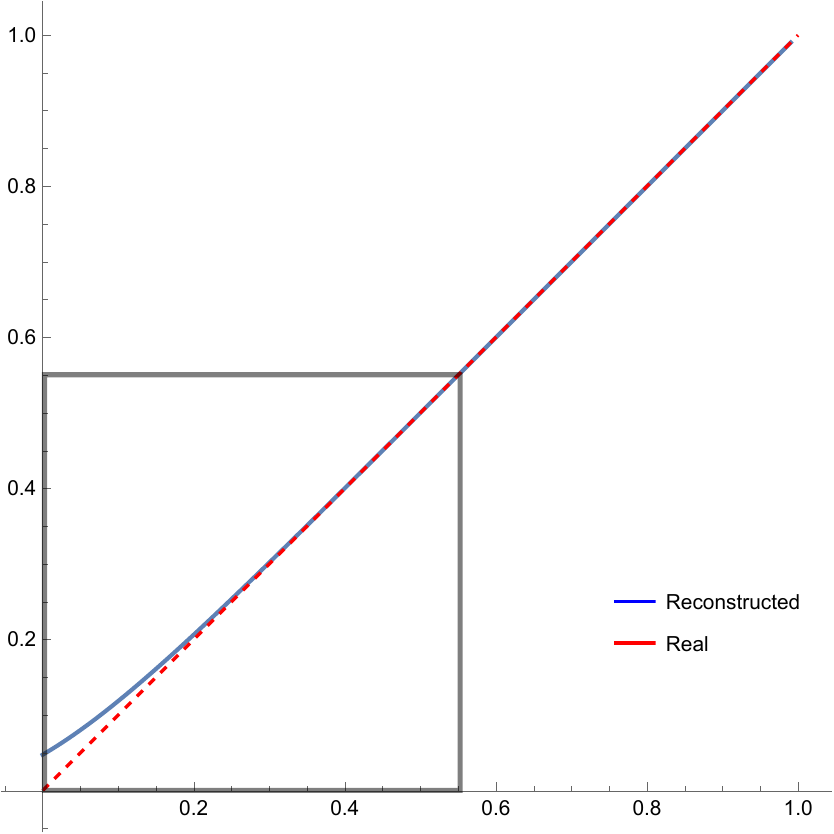}
   %     \caption{Παράδειγμα \ref{ex_teta}.}
      %  \label{fig:sub1teta}
    \end{subfigure}
    \hfill
    % Second subfigure
    \begin{subfigure}[b]{0.45\textwidth}
        \centering
        \includegraphics[width=0.75\textwidth]{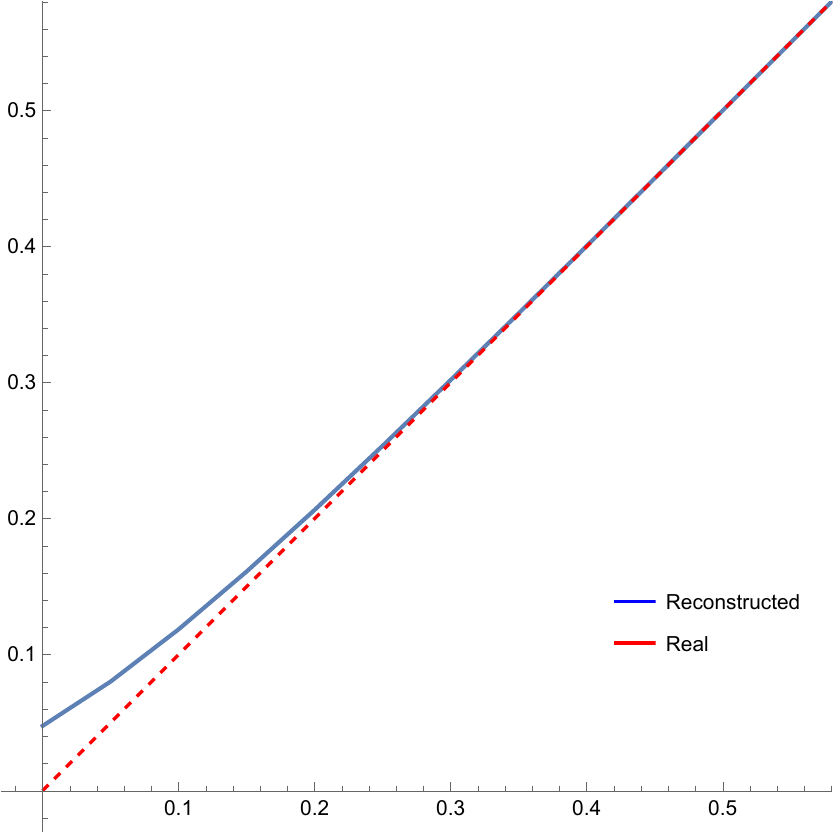}
    %    \caption{Παράδειγμα \ref{ex_trig}.}
     %   \label{fig:sub2teta}
    \end{subfigure}
    \caption{The values of $c:$ real (red) and reconstructed (blue) with $a$ given by  \eqref{eq:solve_a_large}. The right plot is a zoomed-in view of the left.}
    \label{fig:com-c-1}
\end{figure}

Thus, employing the first order approximation yields
\begin{equation}\label{eq:large_a}
I(a) \sim 1 - \frac{4}{\pi} e^{-\frac{\pi^2}{4} a}
\end{equation}
and solving \eqref{eq_final_inverse}, we find
\begin{equation}\label{eq:solve_a_large}
a \approx \frac{4}{\pi^2} \ln\left(\frac{4}{\pi}\frac{1}{1-c} \right),
\end{equation}
which is a sufficient approximation when \(a\) is not small, as illustrated in \autoref{fig:com-c-1}. 

 We define the relative error (RE)
\begin{equation}\label{eq_error}
\mbox{RE} = \frac{|c - \tilde{c} |}{\tilde{c}} \times 100\%,
\end{equation}
for the resulting value of $c$ (normalized data term) using the formula of $a$ and the real $\tilde{c}$ (computed numerically). 

More specifically, the RE using \eqref{eq:solve_a_large} decreases exponentially as \(c \rightarrow 1\), with RE \(= 0.01\%\) at \(c \approx 0.53\). However, for smaller values of \(a\), the approximation may become unreliable. For instance, at  \(c = 0.2\), the RE increases to \(3.2\%\) and at \(c = 0.1\), it rises sharply to \(18.7\%\). In what follows, we aim to achieve small relative error over a larger range of \(c\) either by considering more terms in the Fourier series \eqref{eq:series2} or by employing different approximation schemes.

\begin{remark}\label{ref:bench}
The RE \(= 0.01\%\)
%, which corresponds to \(c \approx 0.53\) and $a \approx 0.4$ 
achieved by the first order Fourier series approximation \eqref{eq:solve_a_large}, will serve as a benchmark for the rest of this paper; we note that this is tighter than the results reported in the literature, see for example \cite{dumm1954, beers65}. Those results correspond (at best) to RE \(= 0.25\%\), obtained at \(c \approx 0.34\) and $a \approx 0.08$.
\end{remark}

\subsubsection{Series inversion}

In this section, we asymptotically solve \eqref{eq_final_inverse} by keeping more terms of the asymptotic series \eqref{eq:series2}, and employing the Lagrange inversion theorem. Setting $a= - \frac{4}{\pi^2} \ln w \Rightarrow w=\exp\left\{ - \frac{\pi^2}{4} a \right\},$  we obtain
\begin{align}\label{eq:series2m}
I(a) = 1 - \frac{4}{\pi}\sum_{n=0}^{\infty} \frac{(-1)^n}{\left( 2n + 1 \right)} w^{ (2n + 1)^2}.
\end{align}
We set $\gamma= \frac{\pi}{4}(1-c)$, then $I(a)=c$ reads
\begin{align}\label{eq:equation2m}
\sum_{n=0}^{\infty} \frac{(-1)^n}{\left( 2n + 1 \right)} w^{ (2n + 1)^2} = \gamma,
\end{align}
 which we invert for $\gamma\to0$, which corresponds to $c\to1$ and $a\to\infty$. 
 
 One should truncate the series, and then use the formula of the Lagrange inversion, known as the Lagrange–B\"urmann formula (for example \cite{WW20} p.129, alternatively \cite{AS65} p.16), which simplifies to the following  form of the solution
\begin{align}
\sum_{n=0}^{N} \frac{(-1)^n}{\left( 2n + 1 \right)} w^{ (2n + 1)^2} = \gamma \ \Rightarrow \ w= \sum_{n=0}^{M} g_n \gamma^{ 8n + 1} + O\left(\gamma\right)^{8M+9},
\end{align}
where $M$ is given such that $8M+1=  (2N + 1)^2$, namely $M=\frac{N(N+1)}{2}$, and $g_0 \equiv 1$.

Then, applying the above truncation and the definition of $w,$ we obtain
\begin{align*}
a&= - \frac{4}{\pi^2} \ln w =   - \frac{4}{\pi^2} \ln \left(\sum_{n=0}^{M} g_n \gamma^{ 8n + 1}\right) =  - \frac{4}{\pi^2} \ln \gamma  - \frac{4}{\pi^2} \ln \left(\sum_{n=0}^{M} g_n \gamma^{ 8n}\right) \notag\\ 
&=   - \frac{4}{\pi^2} \ln \gamma - \frac{4}{\pi^2} \ln \left(1+\sum_{n=1}^{M} g_n \gamma^{ 8n}\right) =    - \frac{4}{\pi^2} \ln \gamma - \frac{4}{\pi^2} \left(\sum_{n=1}^{M} f_n \gamma^{ 8n}\right)  + O\left(\gamma\right)^{8(M+1)}.
\end{align*}
In the last step, we have used the Taylor series expansion. Thus, using the definition of $\gamma$ we find
\begin{equation}\label{sol_fourierN}
a=\frac{4}{\pi^2} \ln \left(\frac{4}{\pi}\frac{1}{1-c}\right) - \frac{4}{\pi^2} \sum_{n=1}^{M} f_n \left[\frac{\pi}{4}(1-c)\right]^{ 8n}  + O\big(1-c\big)^{8(M+1)}.
\end{equation}

\begin{figure}[h!]
    \centering
    % First subfigure
    \begin{subfigure}[b]{0.45\textwidth}
        \centering
        \includegraphics[width=0.75\textwidth]{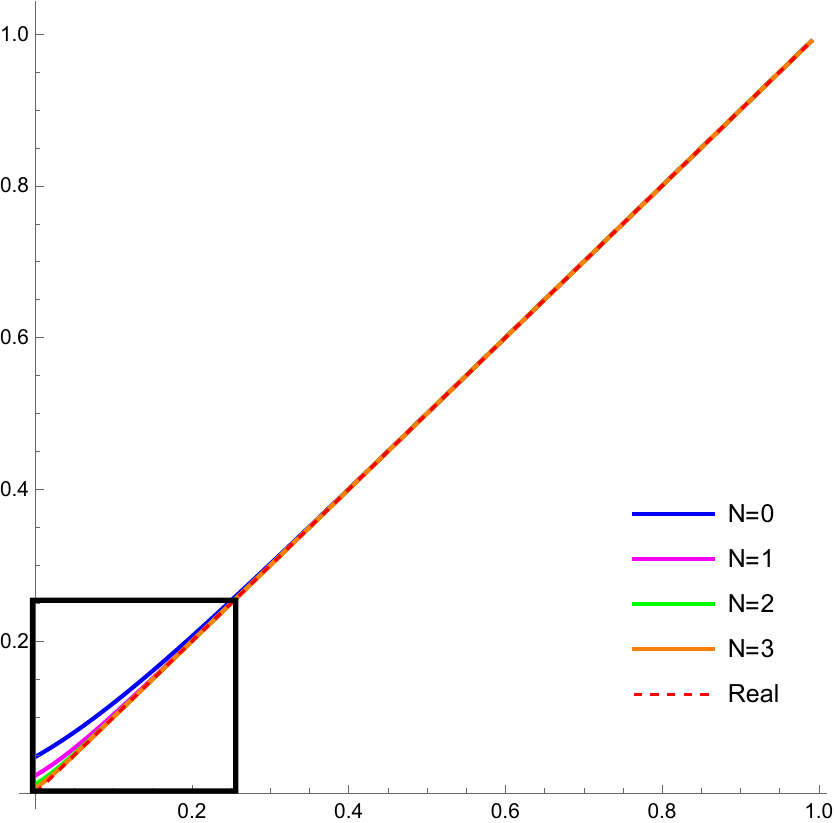}
   %     \caption{Παράδειγμα \ref{ex_teta}.}
      %  \label{fig:sub1teta}
    \end{subfigure}
    \hfill
    % Second subfigure
    \begin{subfigure}[b]{0.45\textwidth}
        \centering
        \includegraphics[width=0.75\textwidth]{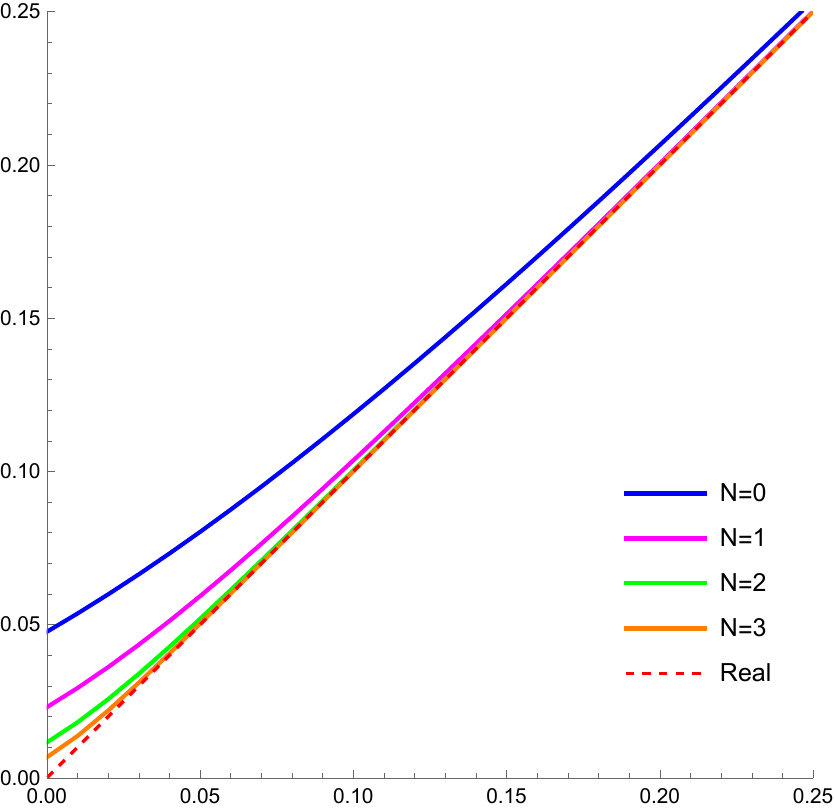}
    %    \caption{Παράδειγμα \ref{ex_trig}.}
     %   \label{fig:sub2teta}
    \end{subfigure}
    \caption{The values of $c:$ real (red) and reconstructed with $a$ given by  \eqref{sol_fourierN} for different values of $N$. The right plot is a zoomed-in view of the left.}
    \label{fig:fourier}
\end{figure}

After some lengthy but straightforward calculations, the constants are found to be
\begin{equation}
g_1=\frac{1}{3}, \ g_2=1, \ g_3= \frac{62}{15}, g_4=\frac{2669}{135}, \  g_5=\frac{13846}{135}, \  g_6=\frac{317783}{567}, 
\end{equation}
and
\begin{equation}\label{eq_f}
f_1=\frac{1}{3}, \ f_2=\frac{17}{18}, \ f_3=\frac{1544}{405}, f_4=\frac{29161}{1620}, \  f_5=\frac{112504}{1215}, \  f_6=\frac{192488308}{382725}.
\end{equation}
Thus, for \(N = 1 \, (M = 1)\), we obtain a second-order approximation; for \(N = 2 \, (M = 3)\), a third-order approximation, and so on. In \autoref{fig:fourier}, we plot the approximated solutions for \(N = 0, 1, 2, 3\), clearly illustrating the improvement of the reconstruction as \(N\) increases. In \autoref{table1}, we summarize the \(c\)-intervals corresponding to the different order approximations, each with a maximum error of approximately \(0.01\%\), in view of \autoref{ref:bench}.

\begin{table}[h!]
\centering
\begin{tabular}{|c | c| c|c|} 
\hline
\(N\) & $M$ & \(c >\) & $a >$ \\
\hline
0 & $-$ & $ 0.531$ & 0.405 \\ 
1 & 1 & $ 0.316$ & 0.251\\ 
2 & 3 & $ 0.172$ & 0.169 \\ 
3 & 6 & $ 0.102 $ & 0.131 \\ 
\hline
\end{tabular}
\caption{Values of \(c\) and $a$ for different \(N\), respectively $M,$ achieving relative error smaller than $0.01\%$ when using \eqref{sol_fourierN}.}\label{table1}
\end{table}

%For numerical purposes
%\begin{equation}
%f_1\approx 0.333, \ f_2\approx 0.944, \ f_3\approx 3.81, f_4\approx 18.0, \  f_5\approx 92.6, \  f_6\approx502.9.
%\end{equation}

\subsubsection{Further approximations}

We present a different approach to improve the estimate of the solution if \( a>0 \) is smaller  than the value $0.405$, reported for \eqref{eq:solve_a_large}.  We set \( a = a^* - \epsilon \), for some small \( \epsilon > 0 \) where $a^\ast$ is given by \eqref{eq:solve_a_large}. Then, we examine the  asymptotic behaviour as \( \epsilon \to 0 \).

From \eqref{eq_integral} we get
\begin{align}\label{eq:small_a_expansion}
I(a) = \frac{i}{\pi} \int_{C}  \frac{e^{-a^* k^{2}}}{k \cos k}  e^{\epsilon k^{2}}dk
     = \sum_{n=0}^{\infty} \frac{i}{\pi} \frac{\epsilon^n}{n!} \int_{C}   \frac{e^{-a^* k^{2}}}{k \cos k} k^{2n} dk = \sum_{n=0}^{\infty} \epsilon^n b_n,  
\end{align}
with the coefficients given by
\begin{equation}\label{eq:bn_def}
b_n = b_n(c) := \frac{i}{\pi} \frac{1}{n!} \int_{C} e^{-a^* k^{2}}  \frac{k^{2n - 1}}{\cos k}  dk = \frac{i}{\pi} \frac{1}{n!} \int_{C} \left(\frac{4}{\pi}\frac{1}{1-c} \right)^{-\frac{4}{\pi^2}k^2}  \frac{k^{2n - 1}}{\cos k}  dk.
\end{equation}

The problem equivalent to \eqref{eq_final_inverse} can be stated as follows: given \( c > 0 \), determine \( \epsilon \) from
\begin{equation}\label{eq:solve_poly}
\sum_{n=0}^{\infty} b_n  \epsilon^n = c. 
\end{equation}

Initially, we consider the second order truncation
\begin{equation}\label{eq:truncation}
I(a) \sim b_0 + b_1  \epsilon + b_2  \epsilon^2,
\end{equation}
to find a closed form solution for $\epsilon>0$:
\begin{equation}\label{eq:quadratic}
b_0 + b_1  \epsilon + b_2  \epsilon^2 = c \quad \Rightarrow \quad  \epsilon = \frac{ -b_1 \pm  \sqrt{b_1^2 - 4 b_2 (b_0 - c)} }{2 b_2},
\end{equation}
with $b_n$ defined in \eqref{eq:bn_def}. We observe numerically that only the ``minus$"$ sign gives $\epsilon \in (0,\, a^\ast)$ and thus,
\begin{equation}\label{a-sol-2nd}
a \sim \frac{4}{\pi^2} \ln\left(\frac{4}{\pi}\frac{1}{1-c} \right)  + \frac{ b_1 + \sqrt{b_1^2 - 4 b_2 (b_0 - c)} }{2 b_2}.
\end{equation}

\begin{figure}[h!]
    \centering
    % First subfigure
    \begin{subfigure}[b]{0.45\textwidth}
        \centering
        \includegraphics[width=0.75\textwidth]{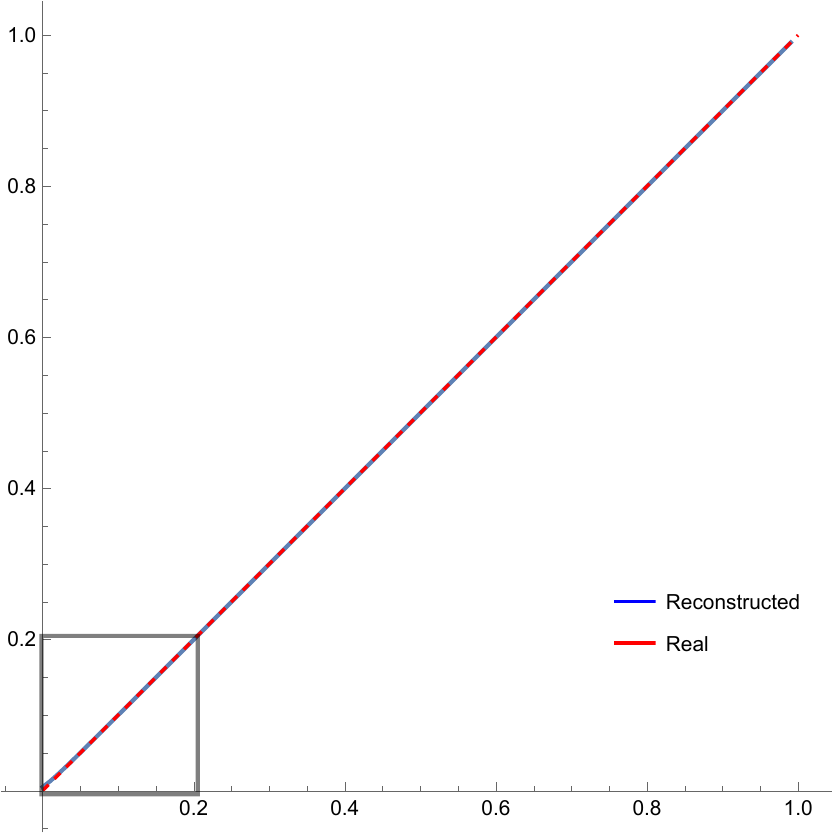}
   %     \caption{Παράδειγμα \ref{ex_teta}.}
      %  \label{fig:sub1teta}
    \end{subfigure}
    \hfill
    % Second subfigure
    \begin{subfigure}[b]{0.45\textwidth}
        \centering
        \includegraphics[width=0.75\textwidth]{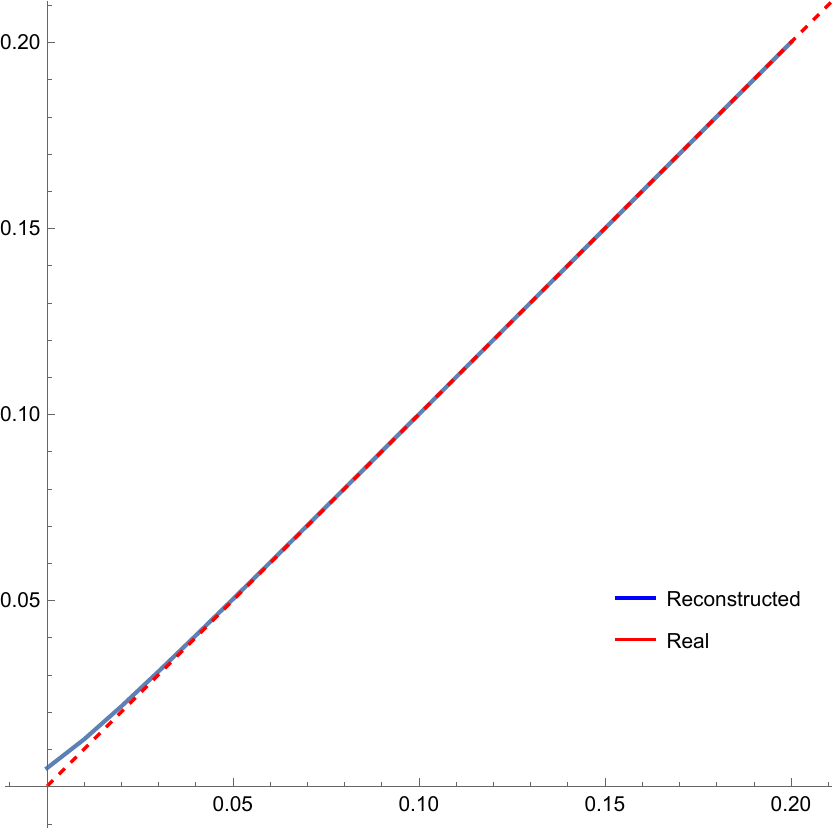}
    %    \caption{Παράδειγμα \ref{ex_trig}.}
     %   \label{fig:sub2teta}
    \end{subfigure}
    \caption{The values of $c:$ real (red) and reconstructed (blue) with $a$ given by  \eqref{a-sol-2nd}. The right plot is a zoomed-in view of the left.}
    \label{fig:quad}
\end{figure}

\autoref{fig:quad} illustrates the approximated $c$ against the real, with $a$ given by \eqref{a-sol-2nd}. It is clear that an almost perfect match is achieved over a larger range of $c$ compared to that obtained for~\eqref{eq:solve_a_large}, and is comparable to the result for~\eqref{sol_fourierN}, when using $N = 3$. In view of  \autoref{ref:bench}, for a relative error below $0.01\%$, we get $c> 0.124$ which corresponds to $a>0.144.$

Of course, one could achieve better results by taking higher order terms in \eqref{eq:solve_poly}. If we consider 5 terms, 
\begin{equation}\label{eq:quartic}
 I(a) \sim \sum_{n=0}^4 b_n \epsilon^n =c, 
\end{equation}
the Ferrari's quartic formula (for example \cite{AS65} p.17) provides an explicit formula for $\epsilon=\epsilon(c)$. Curve fitting improves even more, providing almost perfect reconstructions for  \( c > 0.036 \) and  \( a > 0.089 \), in view of \autoref{ref:bench}.

\subsection{Integral estimation for small $a>0$}\label{sec_small}
We observe that 
$$I(0)=  \frac{i}{\pi} \int_{C}  \frac{1}{k \cos k}  dk =0,$$
by analyticity. However, for $c=0,$ we get $a=0.098$ from \eqref{eq:solve_a_large}, $a = 0.078$ from \eqref{sol_fourierN} for $N=1,$ and $a=0.055$ when considering \eqref{a-sol-2nd}. Thus, it is necessary to obtain approximated solutions which display the property $a \rightarrow 0^+$, as $c \rightarrow0^+.$ Furthermore, employing the analysis of the previous section, we observe that all algebraic higher order terms vanish, indeed, assume $a=\epsilon>0:$  
\begin{align}\label{eq:zero_a_expansion}
I(\epsilon) = \frac{i}{\pi} \int_{C}  \frac{1}{k \cos k}  e^{-\epsilon k^{2}}dk
     = \sum_{n=0}^{\infty} \frac{i}{\pi} \frac{(-\epsilon)^n}{n!} \int_{C}   \frac{1}{k \cos k} k^{2n} dk = \sum_{n=0}^{\infty} \epsilon^n B_n,  
\end{align}
with
\begin{equation*}\label{eq:Bn_def}
B_n = \frac{i}{\pi} \frac{(-1)^n}{n!} \int_{C}   \frac{k^{2n - 1}}{\cos k}  dk =0, \qquad n\in \N.
\end{equation*}
This clearly indicates that the main contribution is exponentially small, hence we accomplish the following procedure
\begin{align*}
I(a)&=I(a)-I(0)=\frac{i}{\pi} \int_{C}  \frac{e^{-a k^{2}}-1}{k \cos k}  dk = \frac{i}{\pi} \int_{C}  \frac{e^{-w^{2}}-1}{w \cos \frac{w}{\sqrt{a}}}  dw, 
\end{align*}
using the change of variables $w=\sqrt{a} \,k$ in the last equality.

Since $C$ lies on the upper complex half plane we apply the expansion
$$
\frac{1}{\cos u}
= \frac{2e^{iu}}{1+e^{2iu}}
= 2e^{iu}\sum_{m=0}^{\infty}(-1)^m e^{2 i m u}
= 2\sum_{m=0}^{\infty}(-1)^m e^{i(2m+1)u},
\quad \text{Im} \, u > 0. 
$$
Thus,
\begin{equation}\label{I-asym-small-a}
\begin{aligned}
I(a)&=\frac{2i}{\pi} \sum_{m=0}^{\infty} (-1)^m \int_{C}   \frac{e^{-w^2}-1}{w}\,e^{i\frac{2m+1}{\sqrt{a}} w}\,dw\\
&=2 \sum_{m=0}^{\infty} (-1)^m \;K\left(\frac{2m+1}{\sqrt{a}}\right)=
2 \sum _{m=0}^{\infty } (-1)^m \text{erfc}\left(\frac{2 m+1}{2 \sqrt{a}}\right),
\end{aligned}
\end{equation}
where we have used that (see \autoref{AppA})
\begin{equation}\label{eq_K}
K(\mu):=\frac{i}{\pi}\int_{-\infty}^{\infty}\frac{e^{-w^2}-1}{w}\,e^{i\mu w}\,dw = \frac{2}{\sqrt{\pi}} \int_{\frac{\mu}{2}}^\infty e^{-z^2} dz  
= : \text{erfc}\left(\frac{\mu }{2}\right),\quad \mu\in\R,
\end{equation}
and  $ \text{erfc}$ denotes the complementary error function

\begin{figure}[h!]
    \centering
    % First subfigure
    \begin{subfigure}[b]{0.45\textwidth}
        \centering
        \includegraphics[width=0.75\textwidth]{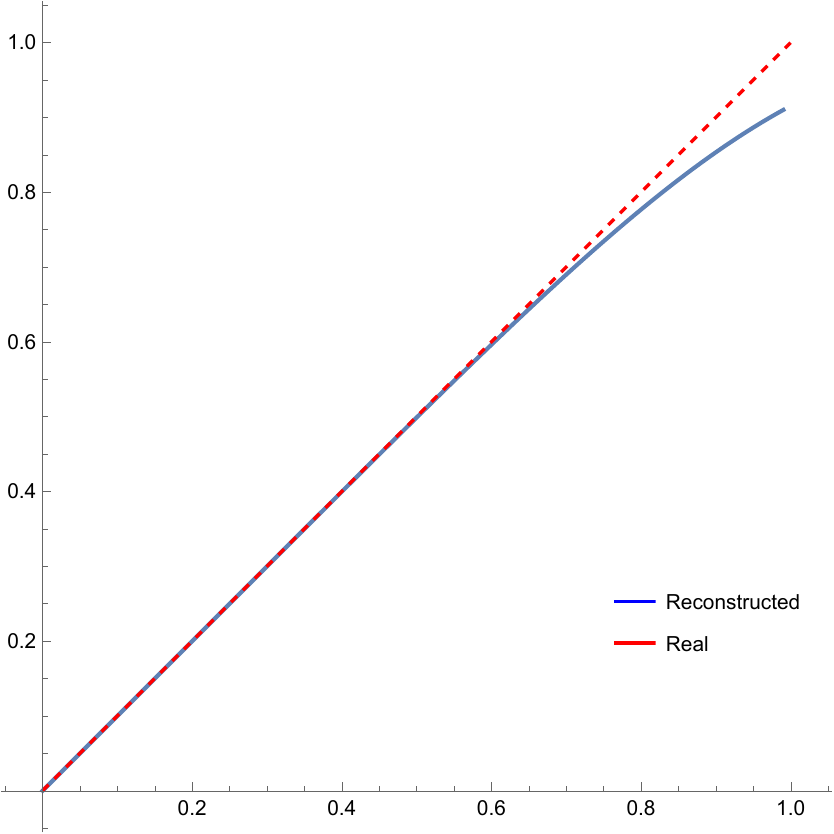}
   %     \caption{Παράδειγμα \ref{ex_teta}.}
      %  \label{fig:sub1teta}
    \end{subfigure}
    \hfill
    % Second subfigure
    \begin{subfigure}[b]{0.45\textwidth}
        \centering
        \includegraphics[width=0.75\textwidth]{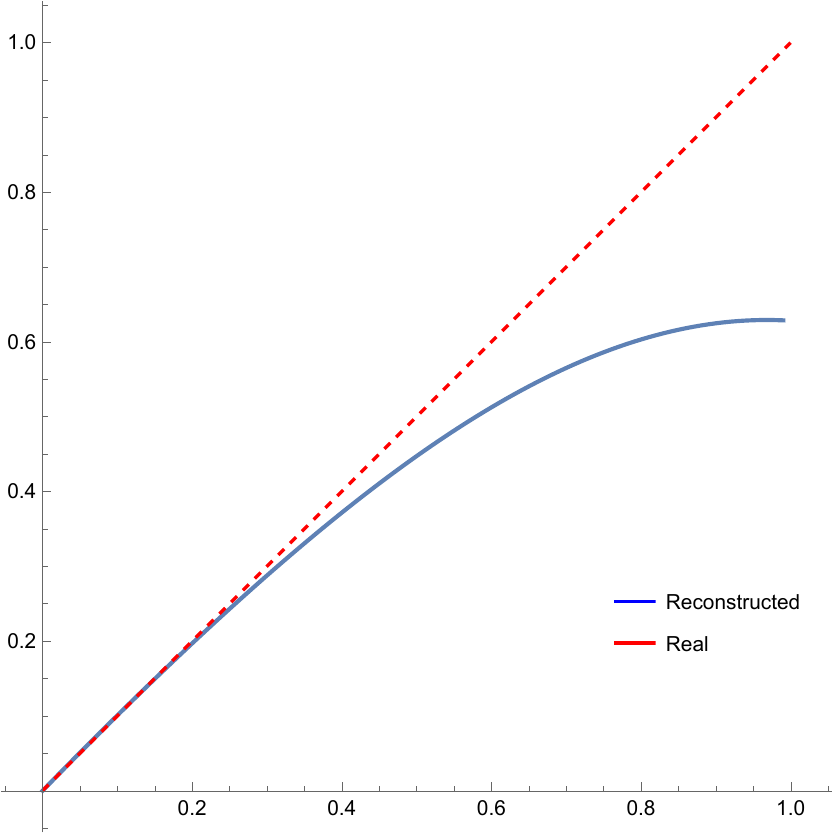}
    %    \caption{Παράδειγμα \ref{ex_trig}.}
     %   \label{fig:sub2teta}
    \end{subfigure}
    \caption{The values of $c:$ real (red) and reconstructed (blue) with $a$ given by \eqref{sol_erf} (left) and by \eqref{a-sol-small-1} (right).}
    \label{fig:erf_log}
\end{figure}

In what follows, we use only the highest order term, namely
\begin{equation}\label{I-asym-small-a-1}
I(a)\sim 2 \, \text{erfc}\left(\frac{1}{2 \sqrt{a}}\right),
\end{equation}
hence, the solution of \eqref{eq_final_inverse} takes the form
\begin{equation}\label{sol_erf}
a\sim\frac{1}{4\, \text{erfc}^{-1}\left(\frac{c}{2}\right)^2},
\end{equation}
where $ \text{erfc}^{-1}$ is the inverse complementary error function. Numerical evaluation  indicates again almost perfect matching for $c< 0.295$, i.e. $a< 0.238$ (see \autoref{ref:bench}). In the left part of \autoref{fig:erf_log} we compare the approximated $c$ with the real, with $a$ given by \eqref{sol_erf}.

In order to get a closed form of \eqref{sol_erf} we employ the asymptotics of $\displaystyle \text{erfc} (z),  \ z\to\infty$, namely
\begin{equation}
\text{erfc} (z) \sim  \frac{e^{-z^2}}{\sqrt{\pi}\,z} 
\sum_{n=0}^\infty (-1)^n \,\frac{(2n-1)!!}{(2z^2)^n},  \quad z\to\infty.
\end{equation}
Keeping the first term, we obtain
$$I(a)\sim \frac{4}{\sqrt{\pi}} \sqrt{a} \, e^{-\frac{1}{4 a}},  \quad a\to0^+.$$
Then,
\[
I(a)=c \ \Leftrightarrow \ \frac{4}{\sqrt{\pi}} \sqrt{a} \, e^{-\frac{1}{4 a}} \sim c, 
\]
which gives
\begin{equation}\label{a-sol-small-1}
a\sim\frac{1}{2 \, W\left(\frac{8}{\pi  c^2}\right)}  \   \Rightarrow  \ a \sim \frac{1}{2 \ln \left[\frac{8}{\pi  c^2} \frac{1}{ \ln \left(\frac{8}{\pi  c^2}\right)}\right]}\,, \quad c\to 0^+,
\end{equation}
where $W$ is the Lambert function. 
%This approximation does not give perfect results as \eqref{sol_erf} and achieves a minimum RE of order $0.09\%$ for $c\approx 0.14.$ 
The property ``$c\to0^+ \Rightarrow a\to0^+$" holds, but the rate of the convergence of this asymptotic expansion is not fully obtained, thus we allow RE up to $3\%,$ which yileds the values of $c<0.258$, which corresponds to $a<0.217$. In the right part of \autoref{fig:erf_log} we plot the approximated $c$ and the real one, with $a$ given by \eqref{a-sol-small-1}.

% Finally, employing the asymptotic expansion of $\frac{1}{4\, \text{erf}^{-1}\left(\frac{2-c}{2}\right)^2}$ as $c\to 0^+$, we obtain 
%\begin{equation}\label{a-sol-small-1}
%a\sim \frac{1}{2 \ln \left[\frac{8}{\pi  c^2} \frac{1}{ \ln \left(\frac{8}{\pi  c^2}\right)}\right]}, \qquad c\to 0^+.
%\end{equation}
%Numerics display almost perfect\footnote{We take as benchmark the values of \cite{Agric-Rich}} match for $c<0.21$, which corresponds to $a<0.19$.

\begin{figure}[h!]
    \centering
    % First subfigure
    \begin{subfigure}[b]{0.45\textwidth}
        \centering
        \includegraphics[width=0.75\textwidth]{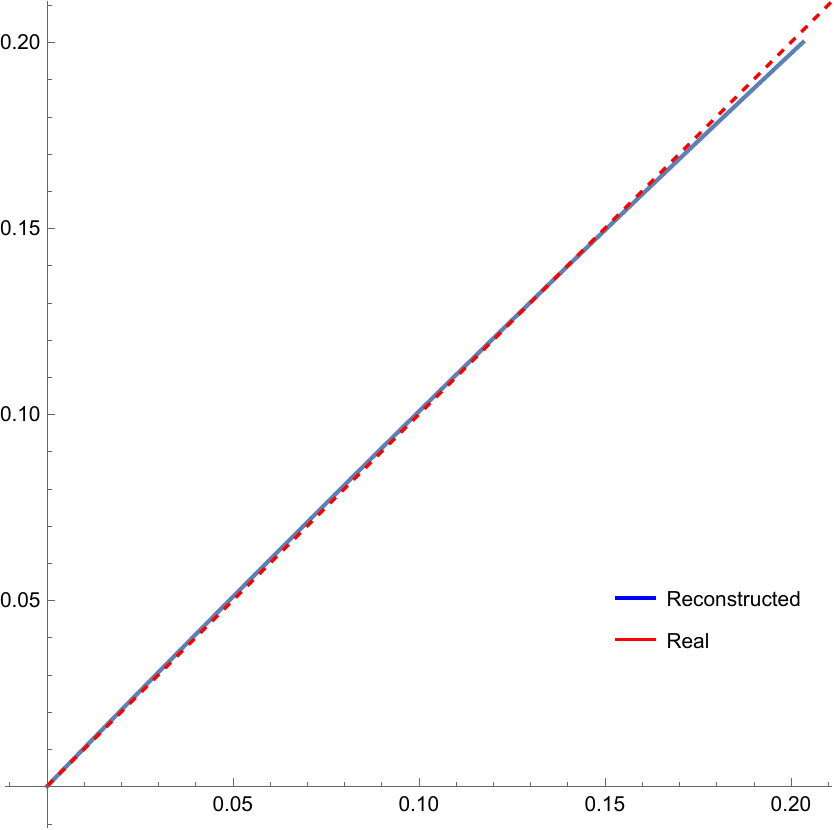}
   %     \caption{Παράδειγμα \ref{ex_teta}.}
      %  \label{fig:sub1teta}
    \end{subfigure}
    \hfill
    % Second subfigure
    \begin{subfigure}[b]{0.45\textwidth}
        \centering
        \includegraphics[width=0.75\textwidth]{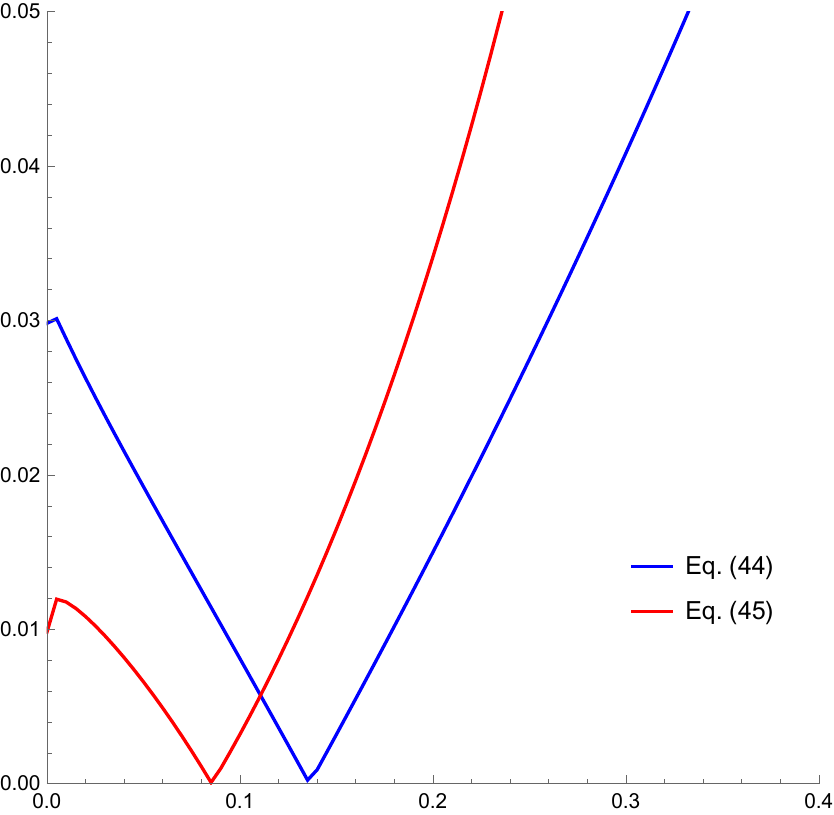}
    %    \caption{Παράδειγμα \ref{ex_trig}.}
     %   \label{fig:sub2teta}
    \end{subfigure}
    \caption{Left: The real (red) and the reconstructed (blue) $c$ values with $a$ given by \eqref{a-sol-small-2}. Right: The relative errors using \eqref{a-sol-small-1} (blue) and \eqref{a-sol-small-2} (red).}
    \label{fig:log_log}
\end{figure}

Finally, by employing a higher-order asymptotic expansion of $\displaystyle \text{erfc} (z),  \ z\to\infty$, we obtain (see \autoref{AppB}) 
\begin{equation}\label{a-sol-small-2}
a\sim \frac{1}{4 P} \left(1+\frac{\ln P}{2 P}+\frac{\ln^2 P-\ln P+2}{4 P^2} \right), \quad P= \ln\frac{2}{c\sqrt{\pi}}\,, \quad c\to 0^+.
\end{equation}
This approximation is more accurate than \eqref{a-sol-small-1} for $c \lesssim 0.11$, with a relative error below $1.2\%$, but the error then increases sharper otherwise (see \autoref{fig:log_log}).

\section{Main Results}\label{sec_main}

In this section we combine the results of \autoref{sec_approx}, namely the asymptotic solutions  $a=a_{as}(c)$  of the integral equation \eqref{eq_final_inverse}, in a way such that the proposed formulas respect the properties 
$$\lim_{a\to 0^+} I(a)=0^+ \qquad \text{and} \qquad \lim_{a\to +\infty} I(a)=1^-,$$
as well as
$$\lim_{c\to 0^+} a(c)=0^+ \qquad \text{and} \qquad \lim_{c\to1^-} a(c)= +\infty,$$
and provide a reliable approximation for the rate of convergence of these limits; in other words, we provide the exact first order term and accurate higher order terms of the relevant asymptotic expansions. In connection to this, the results are presented in terms of the relative error \eqref{eq_error}. 

We combine a formula optimized for large values of $a,$ see \autoref{sec_large}, with another suited for small values, see \autoref{sec_small}, to obtain a uniformly accurate approximation across the full range of $a$.  Specifically, we introduce three expressions: an implicit relation with almost exact agreement, and two explicit approximations whose error bounds remain acceptable relative to their representational complexity.

\begin{figure}[h!]
    \centering

    % First subfigure
    \begin{subfigure}[b]{0.32\textwidth}
        \centering
        \includegraphics[width=\textwidth]{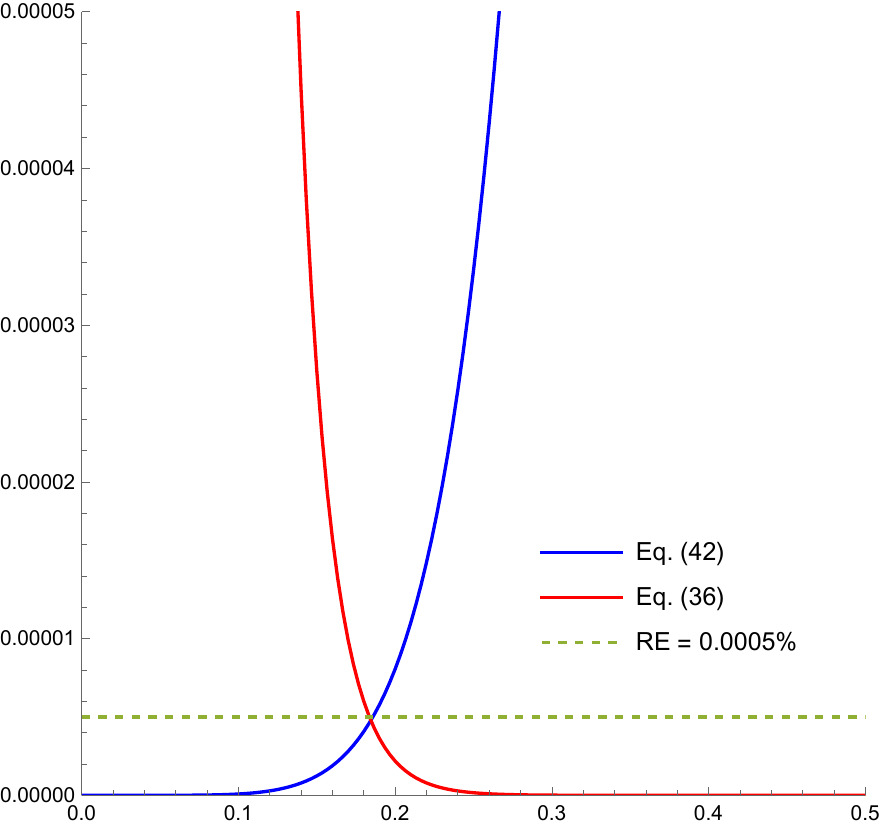}
        \caption{Eq.~\eqref{eq_final1}.}
        \label{fig:com-r-1a}
    \end{subfigure}
    \hfill
    % Second subfigure
    \begin{subfigure}[b]{0.32\textwidth}
        \centering
        \includegraphics[width=\textwidth]{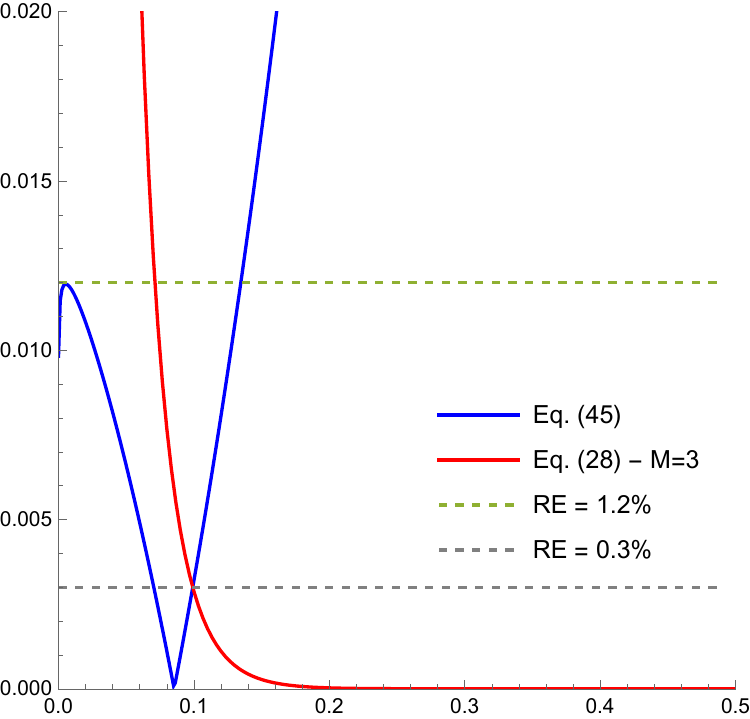}
        \caption{Eq.~\eqref{eq_final2}.}
        \label{fig:com-r-1b}
    \end{subfigure}
    \hfill
    % Third subfigure
    \begin{subfigure}[b]{0.32\textwidth}
        \centering
        \includegraphics[width=\textwidth]{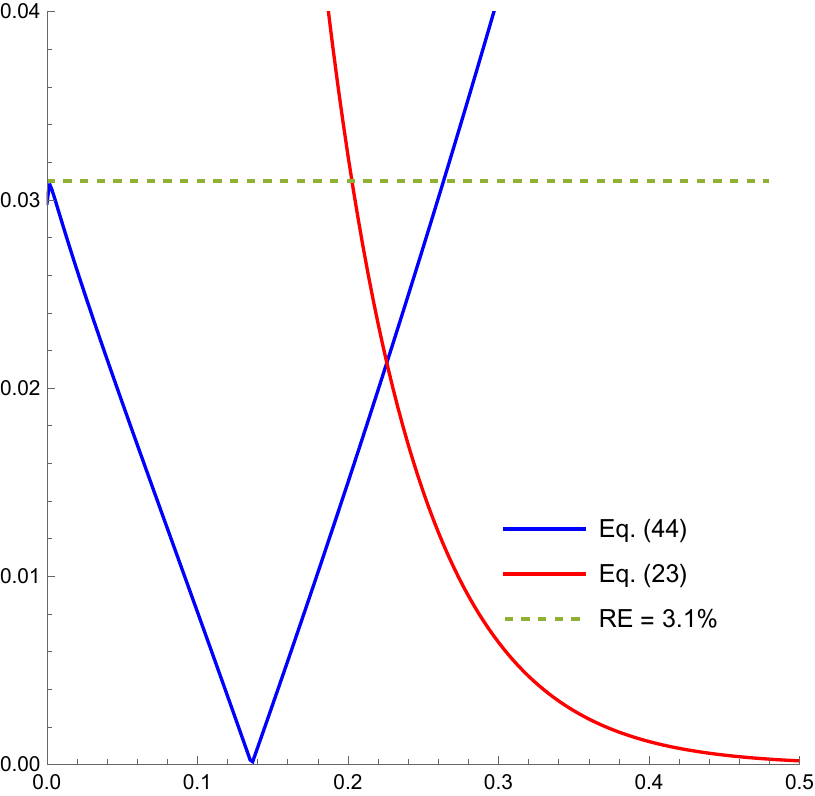}
        \caption{Eq.~\eqref{eq_final3}.}
        \label{fig:com-r-1c}
    \end{subfigure}

    \caption{Relative errors for the different values of $c$. }
    \label{fig:com-r-1}
\end{figure}

\begin{description}
\item[Perfect match]
We combine \eqref{sol_erf} for small $a$ with \eqref{a-sol-2nd} for large values, and we plot the relative error in \autoref{fig:com-r-1a}. The intersection point is at $c \approx 0.18$ with RE$\approx 0.0005\%.$ Thus, we propose:
\begin{equation}\label{eq_final1}
 a = \left.\begin{cases}
			 \frac{1}{4\, \text{erfc}^{-1}\left(\frac{c}{2}\right)^2}, & c\in (0,\,0.18],\\ \vspace*{-3mm}\\
           \frac{4}{\pi^2} \ln\left(\frac{4}{\pi}\frac{1}{1-c} \right)  + \frac{ b_1 + \sqrt{b_1^2 - 4 b_2 (b_0 - c)} }{2 b_2}, & c\in (0.18,1),
		 \end{cases} \right. \quad  \ \mbox{with RE} < 0.0005\%,
\end{equation}
where $ \text{erfc}^{-1}$ denotes the inverse of the complementary error function and $b_n$ are defined in \eqref{eq:bn_def}.

\item[Asymptotic explicit scheme 1] We define the formula:
\begin{equation}\label{eq_final2}
 a = \left.\begin{cases}
			\frac{1}{4 P} \left(1+\frac{\ln P}{2 P}+\frac{\ln^2 P-\ln P+2}{4 P^2} \right), & c\in (0,\,0.1],\\ \vspace*{-3mm}\\
            \frac{4}{\pi^2} \ln \left(\frac{4}{\pi}\frac{1}{1-c}\right) - \frac{4}{\pi^2} \sum_{n=1}^3 f_n \left[\frac{\pi}{4}(1-c)\right]^{ 8n}, & c\in (0.1,\, 1),
		 \end{cases} \right. \quad \ \mbox{with RE} < 1.2\%,
\end{equation}
where $P = \ln \left(\frac{2}{c\sqrt{\pi}}\right)$ and $f_n$ are defined in \eqref{eq_f}. The relative error is presented in \autoref{fig:com-r-1b} where the two curves intersect at $c\approx 0.1$. In fact, the above scheme yields a a maximum RE of $0.3\%$ for $c\in(0.07,1)$. 

\item[Asymptotic explicit scheme 2] This formula offers a closed-form expression that is straightforward to implement while maintaining an acceptable level of relative error, namely:
\begin{equation}\label{eq_final3}
 a = \left.\begin{cases}
			\frac{1}{2 \ln \left[\frac{8}{\pi  c^2} \frac{1}{ \ln \left(\frac{8}{\pi  c^2}\right)}\right]}, & c\in (0,\,0.22],\\ \vspace*{-3mm}\\
            \frac{4}{\pi^2} \ln\left(\frac{4}{\pi}\frac{1}{1-c} \right), & c\in (0.22,\, 1),
		 \end{cases} \right. \quad \ \mbox{with RE} < 3.1\%,
\end{equation}
The intersection point at $c\approx 0.22$ is shown in \autoref{fig:com-r-1c}.

\end{description}

\section{Numerical experiments and physical interpretation}\label{sec_num}

Note that the parameter $a$ depends on $T$,  $L$ and $\mathcal{A}$ (resp. $D_0$). Once $a$ is determined, by keeping two of these parameters fixed, we can address three different inverse problems: two related to the drainage problem (see \autoref{sec_drain}) and one to the infiltration problem (see \autoref{sec_infil}). More precisely, we address the following problems whose solutions are summarized through the solution of \eqref{eq_final_inverse}:

\begin{description}
\item[IP1] Considering \eqref{eq:condition}, find the drain spacing $L,$ given a decrease of the initial height of water table $d+h_0$, to the observed height $H,$ at some time $T.$ We approximate $a$ using any of \eqref{eq_final1}--\eqref{eq_final3} and obtain (see \eqref{def:c+a})
\begin{equation}\label{sol_inv1}
L = \sqrt{\frac{\mathcal{A} T}{a}}.
\end{equation}

\item[IP2] Given $L,$ find the time $T$ needed to reach a specific height $H,$ from \eqref{eq:condition}. The approximated solution is given by any of \eqref{eq_final1}--\eqref{eq_final3} and we get (see \eqref{def:c+a})
\begin{equation}\label{sol_inv2}
T = \frac{a L^2}{\mathcal{A}}.
\end{equation}

\item[IP3] Reconstruct the diffusivity $D_0,$ from one moisture measurement at $x=\frac{L}{2}$ and time $T,$ by solving \eqref{eq:condition2} for $D_0.$ We approximate $a$ by \eqref{eq_final1}--\eqref{eq_final3} and using now \eqref{def:c+a2}, we derive
 \begin{equation}\label{sol_inv3}
D_0 = \frac{a L^2}{4T}.
\end{equation}
\end{description}

For example, using the formula \eqref{eq:solve_a_large}, the solutions of the physically relevant inverse problems are as follows:

\begin{description}
\item[Solution of IP1] Substitute \eqref{eq:solve_a_large} in \eqref{sol_inv1} and using \eqref{eq:condition} we get
\begin{equation}\label{sol_inv1b}
L = \sqrt{\frac{\pi^2 \mathcal{A} T}{4 \ln \left( \frac{4 h_0}{\pi (H-d)} \right) }}.
\end{equation}
\item[Solution of IP2] Similarly, we obtain
\begin{equation}\label{sol_inv2b}
T = \frac{4  L^2 }{\pi^2 \mathcal{A}} \ln \left( \frac{4 h_0}{\pi (H-d)} \right).
\end{equation}
\item[Solution of IP3] Considering now \eqref{eq:solve_a_large} together with \eqref{sol_inv3} and \eqref{eq:condition2} results in
\begin{equation}\label{sol_inv3b}
D_0 = \frac{L^2}{\pi^2 T} \ln \left( \frac{4}{\pi} \frac{\theta_1 - \theta_0}{\theta_1 + \theta_0 - 2 \Theta} \right).
\end{equation}
 
\end{description}

All formulas are well defined, in view of the discussions in \autoref{sec_drain} and \autoref{sec_infil}, respectively. 
The formula \eqref{sol_inv1b} (and, similarly, \eqref{sol_inv2b}) corresponds to the well-known Glover-Dumm equation 
used to compute drain spacing \cite{dumm1954}.  To our knowledge, formula \eqref{sol_inv3b} provides, for the first time, a closed-form expression 
to approximate the diffusivity from a single moisture measurement taken at the midpoint of a bounded, homogeneous soil profile.

\subsection{Drain spacing}

We consider the field data from \cite{kumar} where $h_0 = 1.57$ m. The parameters are summarized in \autoref{table_data1}. Then, the parameter $\mathcal{A}$ is given by $\mathcal{A}= \frac{K B}{S_y},$ where $K$ is the hydraulic conductivity, $B = d+ \frac{h_0}{2},$ and $S_y$ the porosity (dimensionless).

\begin{table}[h!]
\centering
\begin{tabular}{|c|c|c|c|}
\hline
$T$ (days) & $H-d$ (m) & $S_y$ & $K$ (m day$^{-1}$) \\ \hline\hline
1   & 1.38 &  0.060008 & 0.699145  \\
2   & 1.32 &  0.068582 & 0.618233  \\
3   & 1.28 &  0.079471 & 0.577552  \\
4   & 1.24 &  0.083937 & 0.536315  \\
5   & 1.20 &  0.088337 & 0.514509  \\
6   & 1.17 &  0.091103 & 0.474715  \\
7   & 1.13 &  0.091103 & 0.474715  \\
8   & 1.06 &  0.098332 & 0.442264  \\
\hline
\end{tabular}
\caption{Soil drainage properties.}
\label{table_data1}
\end{table}

The authors mention that the actual drain spacing in the field is $2L=50$ m (see \autoref{remark1}), with observed drainage period $T = 2$ days and drop $h_0 - (H-d) = 30$ cm. Nevertheless, the assumptions employed in formulating \eqref{bvp_heat} constrain the model, preventing an exact determination of the actual drain spacing and yielding only an approximate estimate of the optimal spacing $L,$ denoted as ``true$"$ value. 

\begin{table}[h!]
\centering

\begin{tabular}{|c|c|c|c|c|c|c|c|}
\hline
\multirow{2}{*}{$T$ (days)} 
 & \multirow{2}{*}{$c_1$ (Eq. \eqref{def:c+a})} & \multirow{2}{*}{True $2L$} & \multicolumn{4}{|c|}{$2L$ (Eq. \eqref{sol_inv1})}  \\
 \cline{4-7}
& & & Eq. \eqref{eq:solve_a_large} & Eq. \eqref{eq_final1} & Eq.~\eqref{eq_final2} & Eq. \eqref{eq_final3} \\
\hline\hline
1 & 0.12102 & 37.0724 & 36.0371 & 37.0724 & 37.0631 & 37.0426 \\
2 & 0.15924 & 43.0858 & 42.3791 & 43.0858 & 43.0838 & 43.1465 \\
3 & 0.18471 & 47.2058 & 46.6586 & 47.2058 & 47.2050 & 47.3577 \\
4 & 0.21656 & 48.1586 & 47.7961 & 48.1586 & 48.1584 & 48.4420 \\
5 & 0.23567 & 50.3109 & 50.0188 & 50.3109 & 50.3109 & 50.0188 \\
6 & 0.25478 & 51.7023 & 51.4710 & 51.7023 & 51.7023 & 51.4710 \\
7 & 0.28026 & 51.5545 & 51.3923 & 51.5545 & 51.5545 & 51.3923 \\
8 & 0.32484 & 48.4832 & 48.4021 & 48.4832 & 48.4832 & 48.4021 \\
\hline
\end{tabular}
\caption{Comparison of the recovered $L$ values obtained using the different composite approximated solutions of $a,$ introduced in \autoref{sec_main}.}
\label{tab:model_comparison}
\end{table}

Initially, we consider the solution of the first inverse problem, given by \eqref{sol_inv1}. \autoref{tab:model_comparison} compares the reconstructed values of $2L$ obtained from the different composite approximations of $a,$ defined in \autoref{sec_main} with the corresponding true values. The results demonstrate that all three proposed formulations yield accurate reconstructions across the examined range of $T$. Among them, Eq.~\eqref{eq_final1} reproduces the true values with nearly perfect agreement for all cases, confirming its superior accuracy. Eq.~\eqref{eq_final2} also performs well, with deviations smaller than $0.05\%$ even at small $T$ and essentially negligible errors for larger $T$. Furthermore, Eq.~\eqref{eq_final3}, though slightly less precise, still provides reasonable estimates that remain within acceptable error bounds. The Glover--Dumm equation~~\eqref{eq:solve_a_large} underestimates $L$ for smaller $T$, but its accuracy improves noticeably as $T$ increases, consistent with its intended range of validity for large  $c_1.$

\begin{table}[htbp]
\centering
\begin{tabular}{|c|c|c|c|c|c|c|c|}
\hline
\multirow{2}{*}{$2L$ (m)} 
 & \multirow{2}{*}{$c_1$ (Eq.~\eqref{def:c+a})} 
 & \multirow{2}{*}{True $T$} 
 & \multicolumn{4}{|c|}{$T$ (Eq.~\eqref{sol_inv2})}  \\
\cline{4-7}
& & & Eq.~\eqref{eq:solve_a_large} & Eq.~\eqref{eq_final1} & Eq.~\eqref{eq_final2} & Eq.~\eqref{eq_final3} \\
\hline\hline
37 & 0.12102 & 0.9961 & 1.0542 & 0.9961 & 0.9966 & 0.9977 \\
43 & 0.15924 & 1.9920 & 2.0590 & 1.9920 & 1.9922 & 1.9865 \\
47 & 0.18471 & 2.9739 & 3.0441 & 2.9739 & 2.9740 & 2.9549 \\
48 & 0.21656 & 3.9737 & 4.0342 & 3.9737 & 3.9737 & 3.9273 \\
50 & 0.23567 & 4.9384 & 4.9963 & 4.9384 & 4.9384 & 4.9963 \\
52 & 0.25478 & 6.0693 & 6.1240 & 6.0693 & 6.0693 & 6.1240 \\
51 & 0.28026 & 6.8502 & 6.8935 & 6.8502 & 6.8502 & 6.8935 \\
49 & 0.32484 & 8.1714 & 8.1989 & 8.1714 & 8.1714 & 8.1989 \\
\hline
\end{tabular}
\caption{Comparison of the recovered $T$ values obtained using the different composite approximated solutions of $a$, introduced in \autoref{sec_main}.}
\label{tab:model_comparison2}
\end{table}

For completeness, we also present reconstructions of $T,$ given $L,$ namely the second inverse problem with solution given by \eqref{sol_inv2}. In the first column of \autoref{tab:model_comparison2} we set the given $L$ values and in the third column we compute numerically the so-called true $T$ values. The recovered time steps are reported in the last four columns. The Glover--Dumm equation~ \eqref{eq:solve_a_large}, tends to overestimate $T$ for smaller $2L$ but aligns more closely with the true values as $2L$ increases. The results are in accordance with those of the first problem and  demonstrate that the proposed composite approximations deliver reliable reconstructions, with \eqref{eq_final1} offering nearly perfect agreement.

\subsection{Parameter estimation}

We produce simulated data by solving the IBVP \eqref{bc1u} analytically, meaning through evaluating  \eqref{sol_for_infil}. Following \cite{KalMinPal} we consider 9 randomly generated $D_0$ values using
\[
D_0 = 1.2 (1+ \delta) \mbox{ cm}^2/\mbox{h}, 
\] 
with $\delta \in \mathcal{U}(0,1)$ being a uniformly distributed random number in $[0,1].$ We set $\theta_1 = 0.4$ cm$^3$/cm$^3$ and $\theta_0 = 0.05$ cm$^3$/cm$^3$ for a soil profile with length $L=100$ cm.  
% \commentK{Is cm$^3$/cm$^3$ correct? yes}
In the first columns of \autoref{table_infl2}, we summarize the randomized parameter values:  the different time steps $T,$ the data $\Theta$ and the corresponding $c_2$ values. The measurement times are varied to obtain a wider range of $c_2$. 

\begin{table}[h!]
\centering
\scalebox{1}{
\begin{tabular}{|c|c|c|c|c|c|c|c|}
\hline
\multirow{2}{*}{$T$} 
 & \multirow{2}{*}{$\Theta$} 
 & \multirow{2}{*}{$c_2$ (Eq.~\eqref{def:c+a2})} 
 & \multirow{2}{*}{True $D_0$}  
 & \multicolumn{4}{|c|}{$D_0$ (Eq. \eqref{sol_inv3})} \\
 \cline{5-8}
& & & & Eq. \eqref{eq:solve_a_large} & Eq. \eqref{eq_final1} & Eq.~\eqref{eq_final2} & Eq. \eqref{eq_final3} \\
\hline\hline
100  & 0.053097 & 0.017699 & 1.82403 & 2.62849 & 1.82403 & 1.81874 & 1.83671 \\
150  & 0.063641 & 0.077951 & 1.95529 & 2.17990 & 1.95529 & 1.95409 & 1.96542 \\
200  & 0.077174 & 0.155281  & 2.00731 & 2.07868 & 2.00731 & 2.00754 & 2.00270 \\
250  & 0.077733 & 0.158472  & 1.62311 & 1.67828 & 1.62311 & 1.62327 & 1.61871 \\
300  & 0.073859 & 0.136337  & 1.25254 & 1.31088 & 1.25254 & 1.25286 & 1.25251 \\
400  & 0.133329  & 0.476164  & 2.24899 & 2.24969 & 2.24899 & 2.24899 & 2.24969 \\
500  & 0.109254  & 0.338593  & 1.32357 & 1.32721 & 1.32357 & 1.32357 & 1.32721 \\
600  & 0.127599  & 0.443426  & 1.39667 & 1.39742 & 1.39667 & 1.39667 & 1.39742 \\
1000 & 0.197114  & 0.840652  & 2.10569 & 2.10569 & 2.10569 & 2.10569 & 2.10569 \\ \hline
\end{tabular}
}
\caption{Comparison of the reconstructed diffusivity $D_0$ values obtained using the different composite approximated solutions of $a,$ introduced in \autoref{sec_main}. }
\label{table_infl2}
\end{table}

Following the argumentation of the previous subsection, the reconstructed values obtained from \eqref{sol_inv3} for the different formulations of $a_2$ are evaluated using the combined expressions~\eqref{eq_final1}--\eqref{eq_final3} and compared with the true values and \eqref{eq:solve_a_large}. As expected,  eq. \eqref{eq_final1} yields nearly perfect reconstructions;  eq. \eqref{eq_final2} provides highly accurate estimates for $T < 400$ and excellent agreement for $T \ge 400$; eq. \eqref{eq_final3} still produces reasonable results, with the largest deviation observed at $T=100,$ where it overestimates the true $D_0$  by  $0.7\%$;  eq. \eqref{eq:solve_a_large} exhibits its largest error at the earliest time, corresponding to smallest $c_2$ value, which decreases rapidly as $c_2$ increases.

\paragraph*{Author Contributions} All authors contributed equally to this work.

\paragraph*{Funding} The first author acknowledges support by the Sectoral Development Program (SDP 5223471) of the Ministry of Education, Religious Affairs and Sports, through the National Development Program (NDP) 2021-25, grant no. 200/1029.

\paragraph*{Data Availability Statement} The data generated and analyzed during this study are included in this article. Further inquiries can be directed to the corresponding author.

\paragraph*{Conflicts of Interest} The authors declare no conflicts of interest.

\bibliographystyle{siam}
\bibliography{drain_bib}

\appendix
\section*{Appendix}

\section{Proof of \eqref{eq_K} }\label{AppA}
We prove that
$$ K(\mu):=\frac{i}{\pi}\int_{-\infty}^{\infty}\frac{e^{-w^2}-1}{w}\,e^{i\mu w}\,dw = \frac{2}{\sqrt{\pi}} \int_{\frac{\mu}{2}}^\infty e^{-z^2} dz  
= : \text{erfc}\left(\frac{\mu }{2}\right),\qquad \mu\in\R.$$
We differentiate in the distributional sense, w.r.t. to $\mu$, to get
\begin{align*}
\frac{d K}{d\mu} &= \frac{1}{\pi}\int_{-\infty}^{\infty} \ e^{i\mu w}-e^{-w^2+i\mu w}\,dw = \frac{1}{\pi}\int_{-\infty}^{\infty} \ e^{i\mu w}\,dw -\frac{1}{\pi}  \ e^{-\frac{\mu^2}{4}} \int_{-\infty}^{\infty} \ e^{-\left(w^2-i\mu w -\frac{\mu^2}{4}\right)}\,dw \\
&=\delta(\mu)- \frac{1}{\pi}  \ e^{-\frac{\mu^2}{4}}\int_{-\infty}^{\infty} \ e^{-\left(w-i\frac{\mu}{2}\right)^2}\,dw = \delta(\mu)- \frac{1}{\pi}  \ e^{-\frac{\mu^2}{4}}\int_{-\infty-i\frac{\mu}{2}}^{\infty-i\frac{\mu}{2}} \ e^{-z^2}\,dz \\
&=\delta(\mu)- \frac{1}{\sqrt{\pi}}  \ e^{-\frac{\mu^2}{4}} = \delta(\mu) - \frac{d }{d\mu} \text{erf}\left(\frac{\mu }{2}\right).
\end{align*}
Integrating in the distributional sense, w.r.t. to $\mu$, yields the desired result, namely
\begin{align*}
K(\mu)=1 - \text{erf}\left(\frac{\mu }{2}\right)=\text{erfc}\left(\frac{\mu }{2}\right).
\end{align*}

\renewcommand{\theequation}{B.\arabic{equation}}
\setcounter{equation}{0}

\section{Proof of \eqref{a-sol-small-2} }\label{AppB}
We start from the asymptotic equation
\[
c\sim 2\,\text{erfc}\left(\frac{1}{2\sqrt a}\right),\quad a\to0^+ \, .
\]
Letting \(w=\tfrac{1}{4a}\) and using the large-\(w\) expansion of $\text{erfc}$, we obtain
$$\frac{c}{2}\sim\text{erfc}(\sqrt w)\sim \frac{e^{-w}}{\sqrt{\pi w}}\left( 1-\frac{1}{2w}+\frac{3}{4w^2}+\ldots\right).$$
Taking the logarithm of this asymptotic equation yields
\begin{equation}\label{app-main-eq}
w+\frac12\ln w \sim P + r(w),\quad 
P:=\ln\frac{2}{c\sqrt\pi},
\end{equation}
where
\[
r(w):=\ln \left(1-\frac{1}{2w}+\ldots\right) = -\frac{1}{2w}+O(w^{-2}).
\]

We seek for an expansion of the form
\[
w = P-\frac12\ln P + \frac{B\ln P + C}{P} + O\left(\frac{\ln P}{P} \right)^2, \quad \text{for some } \ B,\,C \in\R.
\]
Employing the asymptotics \(\tfrac12\ln w = \tfrac12\ln P -\tfrac14 \tfrac{\ln P }{P}+o(\frac{1}{P})\),  the LHS of \eqref{app-main-eq} reads
\[
w+\frac12\ln w = P + \frac{(B-\tfrac14)\ln P + C}{P} + O\left(\frac{\ln P}{P} \right)^2,
\]
while the RHS of \eqref{app-main-eq} is
\[
P + r(w) = P -\frac{1}{2P} + O\left(\frac{\ln P}{P^2} \right).
\]
Equating coefficients of \(\frac{\ln P}{P}\) and \(\frac{1}{P}\) gives
\[
B=\frac14,\quad C=-\frac12,
\]
resulting in
\[
w = P-\frac12\ln P + \frac{\ln P-2}{4P} + O\left(\frac{\ln P}{P} \right)^2.
\]

Finally, expanding \(a=\tfrac{1}{4w}\) for $P\to\infty$, yields the desired result, namely
\[
a = \frac{1}{4P}\left[1+\frac{\ln P}{2P}+\frac{\ln^2P-\ln P+2}{4P^2}  + O\left(\frac{\ln P}{P} \right)^3\right],
\quad
P=\ln\frac{2}{c\sqrt\pi}, \ c\to0^+.
\]

\end{document}